\documentclass[english,12pt,a4paper]{article}

\usepackage{babel}
\usepackage{epsfig,graphicx,times}   
\usepackage{amssymb}
\newcounter{theoremcounter}
\newcommand{\eop}{\hfill{$\Box$}}
\newcommand{\n}{$\!\!\:\!\!$}
\def\thetheoremcounter{\arabic{theoremcounter}}
\newcommand{\labell}[1]{\label{#1}%
    \ifmmode $$\vspace*{-\baselineskip}\marginpar{#1}%
    \vspace*{-\baselineskip}$$\else\marginpar{#1}\fi}
\newenvironment{satz}%
    {\begin{trivlist}\refstepcounter{theoremcounter}%
    \item[]{{\bf Theorem \thetheoremcounter\ }}\hspace{2mm}\em}%
    {\end{trivlist}}
\newenvironment{lemma}%
    {\begin{trivlist}\refstepcounter{theoremcounter}%
    \item[]{{\bf Lemma \thetheoremcounter\ }}\hspace{2mm}\em}%
    {\end{trivlist}}	 
    {\begin{trivlist}
    \item[]{{\bf Definition
	 \ }}\hspace{2mm}\em}%
    {\end{trivlist}}	 	 
\newenvironment{corollar}%
    {\begin{trivlist}\refstepcounter{theoremcounter}%
    \item[]{{\bf Corollary  \thetheoremcounter\ }}\hspace{2mm}\em}%
    {\end{trivlist}}	 	 	 
    {\begin{trivlist}
    \item[]{\bf Algorithm \thetheoremcounter\ }\hspace{2mm}}%
    {\end{trivlist}}	 
    {\begin{trivlist}
    \item[]{{\sc Proof.}}
    }%
    {\eop\noindent\end{trivlist}}

\begin{document}

\title{Enumeration of generalized polyominoes}
\author{{\sc Matthias Koch} { and }{\sc Sascha Kurz}\\ 
		  Department of Mathematics, University of Bayreuth\\ 
		  matthias.koch, sascha.kurz@uni-bayreuth.de\\
		  D-95440 Bayreuth, Germany}
\maketitle

\vspace*{-7mm}
\noindent
\rule{\textwidth}{0.3 mm}
		 
\begin{abstract}
	\noindent
	As a generalization of polyominoes we consider edge-to-edge connected nonoverlapping unions of regular $k$-gons.    
	For $n\le 4$ we determine formulas for the number $a_k(n)$ of generalized polyominoes consisting of $n$ regular 
	$k$-gons. Additionally we give a table of the numbers $a_k(n)$ for small $k$ and $n$ obtained by computer 
	enumeration. We finish with some open problems for $k$-polyominoes.
\end{abstract}
\noindent
\rule{\textwidth}{0.3 mm}

\section{Introduction}
	A polyomino, in its original definition, is a connected interior-disjoint union of axis-aligned
	unit squares joined edge-to-edge. In other words, it is an edge-connected union of cells in the 
	planar square lattice. For the origin of polyominoes we quote Klarner \cite{handbook_dcg}: 
	\lq\lq Polyominoes have a long 
	
	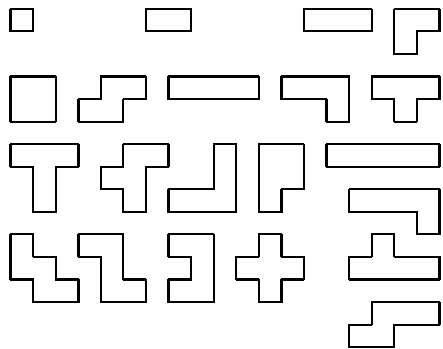
\begin{figure}[ht]
		\begin{center}
			\setlength{\unitlength}{0.30cm}
			\begin{picture}(19,15)
				\put(0,15){\line(1,0){1}}
				\put(0,15){\line(0,-1){1}}
				\put(0,14){\line(1,0){1}}
				\put(1,14){\line(0,1){1}}
				\put(6,14){\line(1,0){2}}
				\put(6,15){\line(1,0){2}}
				\put(6,14){\line(0,1){1}}
				\put(8,14){\line(0,1){1}}
				\put(13,14){\line(1,0){3}}
				\put(13,15){\line(1,0){3}}
				\put(13,14){\line(0,1){1}}
				\put(16,14){\line(0,1){1}}
				\put(17,13){\line(0,1){2}}
				\put(17,13){\line(1,0){1}}
				\put(18,13){\line(0,1){1}}
				\put(18,14){\line(1,0){1}}
				\put(19,14){\line(0,1){1}}
				\put(19,15){\line(-1,0){2}}
				\put(0,12){\line(0,-1){2}}
				\put(0,12){\line(1,0){2}}
				\put(2,10){\line(-1,0){2}}
				\put(2,10){\line(0,1){2}}
				\put(3,11){\line(1,0){1}}
				\put(3,11){\line(0,-1){1}}
				\put(3,10){\line(1,0){2}}
				\put(5,10){\line(0,1){1}}
				\put(5,11){\line(1,0){1}}
				\put(6,11){\line(0,1){1}}
				\put(6,12){\line(-1,0){2}}
				\put(4,12){\line(0,-1){1}}
				\put(7,12){\line(1,0){4}}
				\put(7,11){\line(1,0){4}}
				\put(7,12){\line(0,-1){1}}
				\put(11,12){\line(0,-1){1}}
				\put(12,12){\line(1,0){3}}
				\put(12,12){\line(0,-1){1}}
				\put(12,11){\line(1,0){2}}
				\put(14,11){\line(0,-1){1}}
				\put(14,10){\line(1,0){1}}
				\put(15,10){\line(0,1){2}}
				\put(16,12){\line(1,0){3}}
				\put(16,12){\line(0,-1){1}}
				\put(16,11){\line(1,0){1}}
				\put(17,11){\line(0,-1){1}}
				\put(17,10){\line(1,0){1}}
				\put(18,10){\line(0,1){1}}
				\put(18,11){\line(1,0){1}}
				\put(19,11){\line(0,1){1}}
				\put(0,9){\line(1,0){3}}
				\put(0,9){\line(0,-1){1}}
				\put(0,8){\line(1,0){1}}
				\put(1,8){\line(0,-1){2}}
				\put(1,6){\line(1,0){1}}
				\put(2,6){\line(0,1){2}}
				\put(2,8){\line(1,0){1}}
				\put(3,8){\line(0,1){1}}
				\put(5,9){\line(1,0){2}}
				\put(7,9){\line(0,-1){1}}
				\put(7,8){\line(-1,0){1}}
				\put(6,8){\line(0,-1){2}}
				\put(6,6){\line(-1,0){1}}
				\put(5,6){\line(0,1){1}}
				\put(5,7){\line(-1,0){1}}
				\put(4,7){\line(0,1){1}}
				\put(4,8){\line(1,0){1}}
				\put(5,8){\line(0,1){1}}
				\put(9,9){\line(1,0){1}}
				\put(9,9){\line(0,-1){2}}
				\put(10,9){\line(0,-1){3}}
				\put(10,6){\line(-1,0){3}}
				\put(7,6){\line(0,1){1}}
				\put(7,7){\line(1,0){2}}
				\put(11,9){\line(1,0){2}}
				\put(11,9){\line(0,-1){3}}
				\put(13,9){\line(0,-1){2}}
				\put(13,7){\line(-1,0){1}}
				\put(12,7){\line(0,-1){1}}
				\put(12,6){\line(-1,0){1}}
				\put(19,9){\line(-1,0){5}}
				\put(19,8){\line(-1,0){5}}
				\put(19,9){\line(0,-1){1}}
				\put(14,9){\line(0,-1){1}}
				\put(19,7){\line(-1,0){4}}
				\put(19,7){\line(0,-1){2}}
				\put(19,5){\line(-1,0){1}}
				\put(18,5){\line(0,1){1}}
				\put(18,6){\line(-1,0){3}}
				\put(15,6){\line(0,1){1}}
				\put(0,5){\line(1,0){1}}
				\put(1,5){\line(0,-1){1}}
				\put(1,4){\line(1,0){1}}
				\put(2,4){\line(0,-1){1}}
				\put(2,3){\line(1,0){1}}
				\put(3,3){\line(0,-1){1}}
				\put(3,2){\line(-1,0){2}}
				\put(1,2){\line(0,1){1}}
				\put(1,3){\line(-1,0){1}}
				\put(0,3){\line(0,1){2}}
				\put(3,5){\line(1,0){2}}
				\put(5,5){\line(0,-1){2}}
				\put(5,3){\line(1,0){1}}
				\put(6,3){\line(0,-1){1}}
				\put(6,2){\line(-1,0){2}}
				\put(4,2){\line(0,1){2}}
				\put(4,4){\line(-1,0){1}}
				\put(3,4){\line(0,1){1}}
				\put(7,5){\line(1,0){2}}
				\put(9,5){\line(0,-1){3}}
				\put(9,2){\line(-1,0){2}}
				\put(7,2){\line(0,1){1}}
				\put(7,3){\line(1,0){1}}
				\put(8,3){\line(0,1){1}}
				\put(8,4){\line(-1,0){1}}
				\put(7,4){\line(0,1){1}}
				\put(11,5){\line(1,0){1}}
				\put(12,5){\line(0,-1){1}}
				\put(12,4){\line(1,0){1}}
				\put(13,4){\line(0,-1){1}}
				\put(13,3){\line(-1,0){1}}
				\put(12,3){\line(0,-1){1}}
				\put(12,2){\line(-1,0){1}}
				\put(11,2){\line(0,1){1}}
				\put(11,3){\line(-1,0){1}}
				\put(10,3){\line(0,1){1}}
				\put(10,4){\line(1,0){1}}
				\put(11,4){\line(0,1){1}}
				\put(19,4){\line(0,-1){1}}
				\put(19,4){\line(-1,0){2}}
				\put(19,3){\line(-1,0){4}}
				\put(15,3){\line(0,1){1}}
				\put(15,4){\line(1,0){1}}
				\put(16,4){\line(0,1){1}}
				\put(16,5){\line(1,0){1}}
				\put(17,5){\line(0,-1){1}}
				\put(19,2){\line(-1,0){3}}
				\put(16,2){\line(0,-1){1}}
				\put(16,1){\line(-1,0){1}}
				\put(15,1){\line(0,-1){1}}
				\put(15,0){\line(1,0){2}}
				\put(17,0){\line(0,1){1}}
				\put(17,1){\line(1,0){2}}
				\put(19,1){\line(0,1){1}}
			\end{picture}\\[-3mm]
		\end{center}
		\caption{Polyominoes with at most 5 squares.}
		\label{fig_polyominoes}
	\end{figure}
	
	\noindent
	history, going back to the start
	of the 20th century, but they were popularized in the present era initially
	by Solomon Golomb i.e. \cite{Golomb1954,golomb_first,Golomb1966}, 
	then by Martin Gardner in his \textit{Scientific American} columns.'' 
	At the present time they are widely known by mathematicians, physicists, chemists 
	and have been considered in many different applications, i.e. in the \textit{Ising Model} \cite{Cube}.  
	To give an illustration of polyominoes Figure \ref{fig_polyominoes} depicts the polyominoes consisting 
	of at most 5 unit squares.

	One of the first problems for polyominoes was the determination of there number. Altough there has been some 
	progress, a solution to this problem remains outstanding. In the literature one sometimes speaks also of the 
	cell-growth problem and uses the term animal instead of polyomino.

	Due to its wide area of applications polyominoes were soon generalized to the two 
	other tessellations of the plane, to the eight Archimedean 
	tessellations \cite{archm_tess} and were also considered as unions of $d$-dimensional hypercubes instead of 
	squares. For the known numbers we refer to the \lq\lq Online Encyclopedia of Integer Sequences'' \cite{oeis}. 
	
	\begin{figure}[!h]
		\begin{center}
			\includegraphics{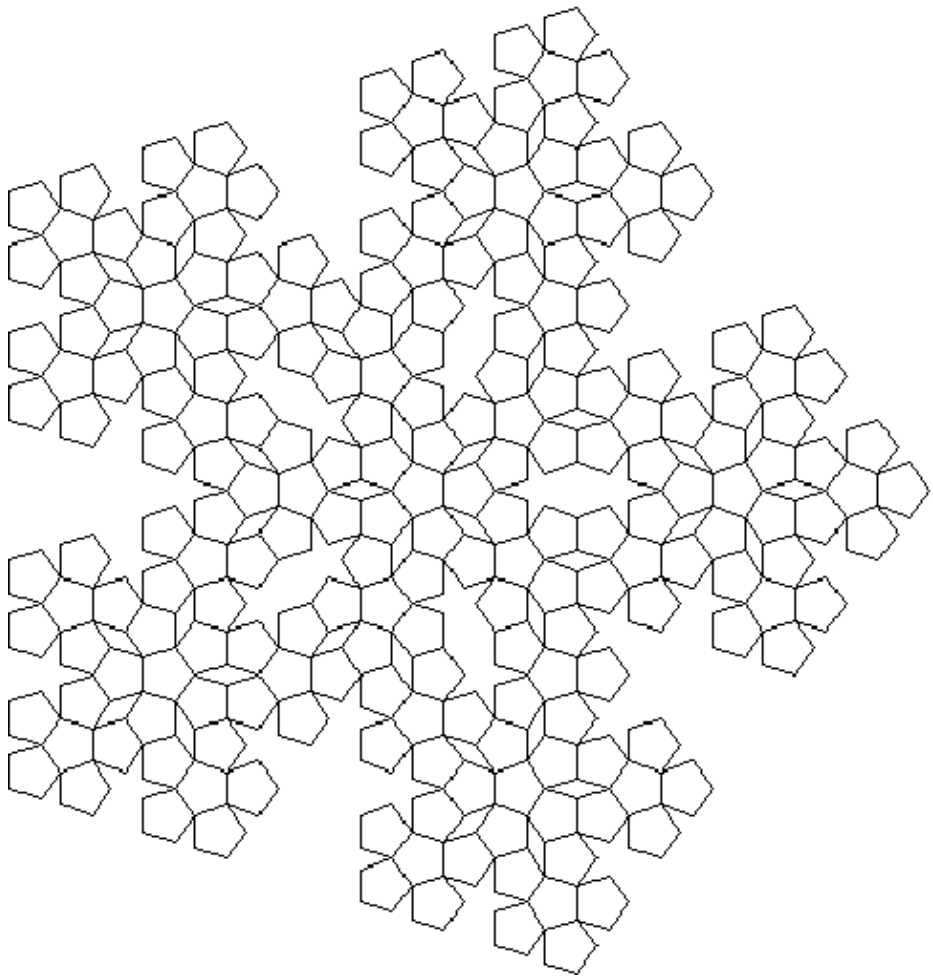}
			\caption{A nice $5$-polyomino.}
			\label{fig_tilling}
		\end{center}		
	\end{figure}
	
	\vspace*{-6mm}
	
	\noindent
	In this article we generalize concept of polyominoes to unions of regular nonoverlapping edge-to-edge connected $k$-gons. 
	For short we call them $k$-polyominoes. An example of a $5$-polyo\-mino, which reminds somewhat to Penrose's famous non-periodic 
	tiling of the plane, is depicted in Figure \ref{fig_tilling}. In the next sections we determine exact formulas for the number 
	$a_k(n)$ of nonisomorphic $k$-polyominoes with $k\le 4$ and give some further values for small parameters $k$ and $n$ obtained 
	by computer enumeration. So far edge-to-edge connected unions of regular $k$-gons were only enumerated if overlapping 
	of the $k$-gons is permitted \cite{0301.05119}. We finish with some open problems for $k$-polyominoes.
		
\section{Formulas for the number of 
$\mathbf{k}$-polyominoes}
	
	By $a_k(n)$ we denote the number of nonisomorphic $k$-polyominoes consisting of $n$ regular $k$-gons as cells
	where $a_k(n)=0$ for $k < 3$. For at most two cells we have $a_k(1)=a_k(2)=1$. If $n\ge 3$ we characterize 
	three edge-to-edge connected cells $\mathcal{C}_1$, 	
	
	\begin{figure}[!ht]
		\begin{center}
			\vspace*{-5mm}
			\includegraphics[width=10cm]{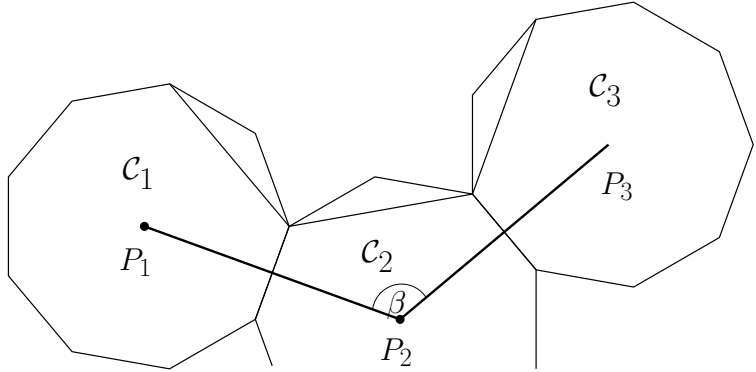}
			\vspace*{-13mm}
			\caption{Angle $\beta=\angle(P_1,P_2,P_3)$ between three neighbored cells.}
			\label{fig_angles_between_3_cells}
		\end{center}		
	\end{figure}
   
	\noindent 
	$\mathcal{C}_2$ and $\mathcal{C}_3$ of a 
	$k$-polyomino, see Figure \ref{fig_angles_between_3_cells}, by the angle $\beta=\angle(P_1,P_2,P_3)$ between the 
	centers of the cells. Since these angles are multiples of $\frac{2\pi}{k}$ we call the minimum 
	$$\min\left( \angle(P_1,P_2,P_3)\frac{k}{2\pi}, \left(2\pi-\angle(P_1,P_2,P_3)\right)
	\frac{k}{2\pi}\right)$$ the discrete angle between $\mathcal{C}_1$, $\mathcal{C}_2$, and $\mathcal{C}_3$ and denote 
	it by $\delta(\mathcal{C}_1,\mathcal{C}_2,\mathcal{C}_3)$.

\begin{figure}[ht]
		\begin{center}
			\vspace*{-10mm}
			\includegraphics{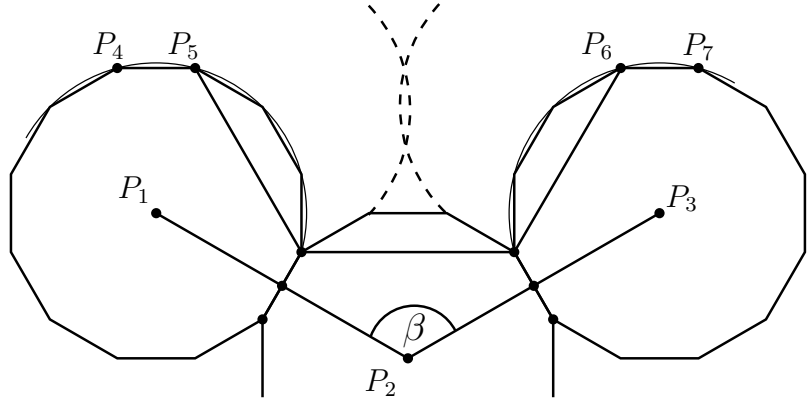}
			\vspace*{-16mm}
			\caption{Nonoverlapping 12-gons.}
			\label{Proof_lemma_one}
		\end{center}	
	\end{figure}

	\begin{lemma}
		\label{lemma_overlapping}
		Two $k$-gons $\mathcal{C}_1$ and $\mathcal{C}_3$ joined via an edge to a $k$-gon $\mathcal{C}_2$ are nonoverlapping 
		if and only if $\delta(\mathcal{C}_1,\mathcal{C}_2,\mathcal{C}_3)\ge\Big\lfloor\frac{k+5}{6}\Big\rfloor$. The three $k$-gons 
		are neighbored pairwise if and only if $k\equiv 0\,\mbox{mod}\,6$.
	\end{lemma}
	\textbf{Proof.}
	We consider Figure \ref{fig_angles_between_3_cells} and set 
	$\beta=\delta(\mathcal{C}_1,\mathcal{C}_2,\mathcal{C}_3)\frac{2\pi}{k}$. If the cells $\mathcal{C}_1$ and $\mathcal{C}_3$ 
	are non-overlapping we have $\overline{P_1P_3}\ge \overline{P_1P_2}$ because the lengths of the lines $\overline{P_1P_2}$ 
	and $\overline{P_2P_3}$ are equal. Thus $\beta\ge \frac{2\pi}{6}$ and 
	$\delta(\mathcal{C}_1,\mathcal{C}_2,\mathcal{C}_3)\ge\Big\lfloor\frac{k+5}{6}\Big\rfloor$ is  necessary. Now we consider 
	the circumcircles of the cells $\mathcal{C}_1$ and $\mathcal{C}_3$, see Figure \ref{Proof_lemma_one}. Due to 
	$\beta\ge\frac{2\pi}{6}$ only the circlesegments between points $P_4$, $P_5$ and $P_6$, $P_7$ may intersect. The last step 
	is to check that the corresponding lines $\overline{P_4P_5}$ and $\overline{P_6P_7}$ do not intersect and they touch each 
	other if and only if $k\equiv 0\,\mbox{mod}\,6$.\hfill{$\square$}
	
	\begin{corollar}
		\label{cor_neighbors}
		The number of neighbors of a cell in a $k$-polyomino is at most $$\min\left(k,\frac{k}{\Big\lfloor\frac{k+5}{6}
		\Big\rfloor}\right)\le 6\,.$$
	\end{corollar}
	
	\noindent
	With the aid of Lemma \ref{lemma_overlapping} we are able to determine the number $a_k(3)$ of $k$-polyominoes consisting 
	of $3$ cells.
	
	\begin{satz}
		\label{theorem_a_k_3}
		$$a_k(3)=\left\lfloor\frac{k}{2}\right\rfloor-\left\lfloor\frac{k+5}{6}\right\rfloor+1\quad\mbox{for}\,\, k\ge 3\,.$$
	\end{satz}
	\textbf{Proof.} It suffices to determine the possible values for $\delta(\mathcal{C}_1,\mathcal{C}_2,\mathcal{C}_3)$. Due 
	to Lemma \ref{lemma_overlapping} we have $\delta(\mathcal{C}_1,\mathcal{C}_2,\mathcal{C}_3)\ge 
	\left\lfloor\frac{k+5}{6}\right\rfloor$ and due to to symmetry considerations we have 
	$\delta(\mathcal{C}_1,\mathcal{C}_2,\mathcal{C}_3)\le \left\lfloor\frac{k}{2}\right\rfloor$.\hfill{$\square$}
	
	In order to determine the number of $k$-polyominoes with more than $3$ cells we describe the classes of $k$-polyominoes
	by graphs. We represent each $k$-gon by a vertex and join two vertices exactly if they are connected via an edge.
	
	\begin{figure}[ht]
		\begin{center}
			\setlength{\unitlength}{0.6cm}
			\begin{picture}(14.5,2)
				\put(0,1){\line(1,0){1}}
				\put(0,1){\circle*{0.3}}
				\put(1,1){\circle*{0.3}}
				\put(1,1){\line(1,1){0.74}}
				\put(1,1){\line(1,-1){0.74}}
				\put(1.74,1.74){\circle*{0.3}}
				\put(1.74,0.26){\circle*{0.3}}
				\put(3,1){\line(1,0){1}}
				\put(3,1){\circle*{0.3}}
				\put(4,1){\circle*{0.3}}
				\put(4,1){\line(1,1){0.74}}
				\put(4,1){\line(1,-1){0.74}}
				\put(4.74,1.74){\circle*{0.3}}
				\put(4.74,0.26){\circle*{0.3}}
				\qbezier(3,1)(3.87,1.37)(4.74,1.74)
				\put(6,1){\line(1,0){1}}
				\put(6,1){\circle*{0.3}}
				\put(7,1){\circle*{0.3}}
				\put(7,1){\line(1,1){0.74}}
				\put(7,1){\line(1,-1){0.74}}
				\put(7.74,1.74){\circle*{0.3}}
				\put(7.74,0.26){\circle*{0.3}}
				\qbezier(6,1)(6.87,1.37)(7.74,1.74)
				\qbezier(7.74,1.74)(7.74,1)(7.74,0.26)
				\put(9,1){\line(1,0){3}}
				\put(9,1){\circle*{0.3}}
				\put(10,1){\circle*{0.3}}
				\put(11,1){\circle*{0.3}}
				\put(12,1){\circle*{0.3}}
				\put(13.5,0.5){\line(1,0){1}}
				\put(13.5,1.5){\line(1,0){1}}
				\put(13.5,0.5){\line(0,1){1}}
				\put(14.5,0.5){\line(0,1){1}}
				\put(13.5,0.5){\circle*{0.3}}
				\put(13.5,1.5){\circle*{0.3}}
				\put(14.5,0.5){\circle*{0.3}}
				\put(14.5,1.5){\circle*{0.3}}
			\end{picture}
			\caption{The possible graphs of $k$-polyominoes with $4$ vertices.}
			\label{fig_graphs}	
		\end{center}	
	\end{figure}
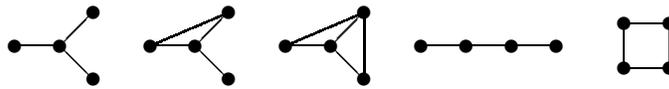
	
	\begin{lemma}
	   \label{lemma_count_4_1}
		The number of $k$-polyominoes with a graph isomorphic to one of the first three ones in Figure \ref{fig_graphs} is given by
		$$
			\left\lfloor\frac{\left(k-3\left\lfloor\frac{k+5}{6}\right\rfloor\right)^2+
			6\left(k-3\left\lfloor\frac{k+5}{6}\right\rfloor\right)+12}{12}\right\rfloor.
		$$
	\end{lemma}
	\textbf{Proof.} We denote the cell corresponding to the unique vertex of degree $3$ in the graph by $\mathcal{C}_0$ 
	and the three other cells by $\mathcal{C}_1$, $\mathcal{C}_2$, and $\mathcal{C}_3$. With
	$\delta_1=\delta(\mathcal{C}_1,\mathcal{C}_0,\mathcal{C}_2)-\left\lfloor\frac{k+5}{6}\right\rfloor$,  
	$\delta_2=\delta(\mathcal{C}_2,\mathcal{C}_0,\mathcal{C}_3)-\left\lfloor\frac{k-1}{6}\right\rfloor$, and
	$\delta_3=\delta(\mathcal{C}_3,\mathcal{C}_0,\mathcal{C}_1)-\left\lfloor\frac{k-1}{6}\right\rfloor$ we set
	$m=\delta_1+\delta_2+\delta_3=k-3\left\lfloor\frac{k+5}{6}\right\rfloor$. Because the $k$-polyominoes with a 
	graph isomorphic to one of the first three ones in Figure \ref{fig_graphs} are uniquely described by 
	$\delta_1,\delta_2,\delta_3$, due to Lemma \ref{lemma_overlapping} and due to symmetry their number
	equals the number of partitions of $m$ into at most three parts. This number is the coefficient 
	of $x^m$ in the Taylor series of $\frac{1}{(1-x)(1-x^2)(1-x^3)}$ in $x=0$ and can be expressed as
	$\left\lfloor\frac{m^2+6m+12}{12}\right\rfloor$. \hfill{$\square$}
	
	\begin{figure}[ht]
		\begin{center}
			\setlength{\unitlength}{0.6cm}
			\begin{picture}(7.5,2)
				\put(0,1){\line(1,-1){1}}
				\put(1,0){\line(1,0){1.5}}
				\put(2.5,0){\line(1,1){1}}
				\put(0,1){\circle*{0.3}}
				\put(1,0){\circle*{0.3}}
				\put(2.5,0){\circle*{0.3}}
				\put(3.5,1){\circle*{0.3}}
				\qbezier(0.6,0.4)(1.25,0.7)(1.57,0)
				\qbezier(1.93,0)(2.15,0.7)(2.9,0.4)
				\put(4.5,0){\line(1,1){1}}
				\put(5.5,1){\line(1,-1){1}}
				\put(6.5,0){\line(1,1){1}}
				\put(4.5,0){\circle*{0.3}}
				\put(5.5,1){\circle*{0.3}}
				\put(6.5,0){\circle*{0.3}}
				\put(7.5,1){\circle*{0.3}}
				\qbezier(5.1,0.6)(5.5,0.3)(5.9,0.6)
				\qbezier(6.1,0.4)(6.5,0.7)(6.9,0.4)
			\end{picture}
			\caption{Paths of lengths $3$ representing chains of four neighbored cells.}
			\label{fig_paths}	
		\end{center}	
	\end{figure}
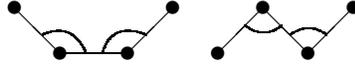
	
	\noindent
	In Lemma \ref{lemma_overlapping} we have given a condition for a chain of three neighbored cells avoiding an overlapping. For
	a chain of four neighbored cells we have to consider the two cases of Figure \ref{fig_paths}. In the second case the two 
	vertices of degree one are not able to overlap so we need a lemma in the spirit of Lemma \ref{lemma_overlapping} only for the 
	first case.
	
	\begin{lemma}
		\label{lemma_overlap_2}
		Four $k$-gons $\mathcal{C}_1$, $\mathcal{C}_2$, $\mathcal{C}_3$, and $\mathcal{C}_4$ arranged as in the first case 
		of Figure \ref{fig_paths} are nonoverlapping if and only if Lemma \ref{lemma_overlapping} is fulfilled for the two subchains of
		length $3$ and
		$$
			\delta(\mathcal{C}_1,\mathcal{C}_2,\mathcal{C}_3)+\delta(\mathcal{C}_2,\mathcal{C}_3,\mathcal{C}_4)\ge
			\left\lfloor\frac{k+1}{2}\right\rfloor\,.
		$$
		The chain is indeed a $4$-cycle if and only if 
		$$
			\delta(\mathcal{C}_1,\mathcal{C}_2,\mathcal{C}_3)+\delta(\mathcal{C}_2,\mathcal{C}_3,\mathcal{C}_4)
			=\frac{k}{2}\,.
		$$ 
	\end{lemma}
	\textbf{Proof.} We start with the second statement and consider the quadrangle of the centers of the $4$ 
	cells. Because the angle sum of a quadrangle is $2\pi$ we have  
	$$
		\delta(\mathcal{C}_1,\mathcal{C}_2,\mathcal{C}_3)+
		\delta(\mathcal{C}_2,\mathcal{C}_3,\mathcal{C}_4)+d(\mathcal{C}_3,\mathcal{C}_4,\mathcal{C}_1)+
		\delta(\mathcal{C}_4,\mathcal{C}_1,\mathcal{C}_2)=k\,.
	$$ 
	Due to the fact that the side lengths of the quadrangle are equal  we have
	$$
		\delta(\mathcal{C}_1,\mathcal{C}_2,\mathcal{C}_3)+
		\delta(\mathcal{C}_2,\mathcal{C}_3,\mathcal{C}_4)=
		\delta(\mathcal{C}_3,\mathcal{C}_4,\mathcal{C}_1)+
		\delta(\mathcal{C}_4,\mathcal{C}_1,\mathcal{C}_2)
	$$
	which is equivalent to the statement.
	
	Thus $\delta(\mathcal{C}_1,\mathcal{C}_2,\mathcal{C}_3)+\delta(\mathcal{C}_2,\mathcal{C}_3,\mathcal{C}_4)\ge
	\left\lfloor\frac{k+1}{2}\right\rfloor$ is a necessary condition. Similar to the proof of Lemma \ref{lemma_overlapping} 
	we consider the circumcircles of the cells $\mathcal{C}_1$, $\mathcal{C}_4$ and check that the cells do not 
	intersect.\hfill{$\square$}
	
	\begin{lemma}
		\label{lemma_count_4_2}
		For $k\ge 3$ the number of $k$-polyominoes with a graph isomorphic to one of the last two ones in Figure \ref{fig_graphs} is
		given by\\[-3mm]
		\begin{eqnarray*}
			\frac{5k^2+4k}{48}\,\,\,\mbox{for}\,\,\, k\equiv 0\,\mbox{mod}\,12,&\quad&
			\frac{5k^2+6k-11}{48}\,\,\,\mbox{for}\,\,\, k\equiv 1\,\mbox{mod}\,12,\\
			\frac{5k^2+12k+4}{48}\,\,\,\mbox{for}\,\,\, k\equiv 2\,\mbox{mod}\,12,&\quad&
			\frac{5k^2+14k+9}{48}\,\,\,\mbox{for}\,\,\, k\equiv 3\,\mbox{mod}\,12,\\
			\frac{5k^2+20k+32}{48}\,\,\,\mbox{for}\,\,\, k\equiv 4\,\mbox{mod}\,12,&\quad&
			\frac{5k^2+22k+5}{48}\,\,\,\mbox{for}\,\,\, k\equiv 5\,\mbox{mod}\,12,\\
			\frac{5k^2+4k-12}{48}\,\,\,\mbox{for}\,\,\, k\equiv 6\,\mbox{mod}\,12,&\quad&
			\frac{5k^2+6k+1}{48}\,\,\,\mbox{for}\,\,\, k\equiv 7\,\mbox{mod}\,12,\\
			\frac{5k^2+12k+16}{48}\,\,\,\mbox{for}\,\,\, k\equiv 8\,\mbox{mod}\,12,&\quad&
			\frac{5k^2+14k-3}{48}\,\,\,\mbox{for}\,\,\, k\equiv 9\,\mbox{mod}\,12,\\
			\frac{5k^2+20k+20}{48}\,\,\,\mbox{for}\,\,\, k\equiv 10\,\mbox{mod}\,12,&\quad&
			\frac{5k^2+22k+17}{48}\,\,\,\mbox{for}\,\,\, k\equiv 11\,\mbox{mod}\,12\,.
		\end{eqnarray*}
	\end{lemma}
	\textbf{Proof.} Because each of the last two graphs in Figure \ref{fig_graphs} contains a path of length $3$ as a subgraph
	we consider the two cases of Figure \ref{fig_paths}. We denote the two interesting discrete angles by $\delta_1$ and 
	$\delta_2$. Due to symmetry we may assume $\delta_1\le \delta_2$ and because the graphs do not contain a triangle we have 
	$\delta_2\ge\delta_1\ge\left\lfloor\frac{k+6}{6}\right\rfloor$ due to Lemma \ref{lemma_overlapping}.
	From the definition of the discrete angle we have $\delta_1\le\delta_2\le\left\lfloor\frac{k}{2}\right\rfloor$. To avoid 
	double counting we assume $\left\lfloor\frac{k+6}{6}\right\rfloor\le \delta_1\le\delta_2\le\left\lfloor\frac{k-1}{2}\right\rfloor$ 
	in the second case, so that we get a number of 
	$$
		{\left\lfloor\frac{k-1}{2}\right\rfloor-\left\lfloor\frac{k+6}{6}\right\rfloor+2\choose 2}
	$$
	$k$-polyominoes. With Lemma \ref{lemma_overlap_2} and a look at the possible symmetries the number of $k$-polyo\-minoes 
	in the first case is given by 
	$$
		\sum_{\delta_1=\left\lfloor\frac{k+6}{6}\right\rfloor}^{\left\lfloor\frac{k}{2}\right\rfloor}
		\sum_{\delta_2=\max(\delta_1,\left\lfloor\frac{k+1}{2}\right\rfloor-\delta_1)}^{\left\lfloor\frac{k}{2}\right\rfloor}1\,.
	$$
	A little calculation yields the proposed formulas. \hfill{$\square$}
	
	\begin{satz}
		\label{theorem_a_k_4}
		For $k\ge 3$ we have
		$$
			a_k(4)=\left\{
			\begin{array}{cccccc}
 				\frac{3k^2+8k+24}{24}  & \mbox{for} & k\equiv 0\,\,\mbox{mod}\,\,12,&
				\frac{3k^2+4k-7}{24}   & \mbox{for} & k\equiv 1\,\,\mbox{mod}\,\,12, \\ 
 				\frac{3k^2+8k-4}{24}   & \mbox{for} & k\equiv 2\,\,\mbox{mod}\,\,12, & 
 				\frac{3k^2+10k+15}{24} & \mbox{for} & k\equiv 3\,\,\mbox{mod}\,\,12, \\ 
			 	\frac{3k^2+14k+16}{24} & \mbox{for} & k\equiv 4\,\,\mbox{mod}\,\,12, & 
 				\frac{3k^2+16k+13}{24} & \mbox{for} & k\equiv 5\,\,\mbox{mod}\,\,12, \\ 
 				\frac{3k^2+8k+12}{24}  & \mbox{for} & k\equiv 6\,\,\mbox{mod}\,\,12, & 
 				\frac{3k^2+4k-7}{24}   & \mbox{for} & k\equiv 7\,\,\mbox{mod}\,\,12, \\ 
 				\frac{3k^2+8k+8}{24}   & \mbox{for} & k\equiv 8\,\,\mbox{mod}\,\,12, & 
 				\frac{3k^2+10k+3}{24}  & \mbox{for} & k\equiv 9\,\,\mbox{mod}\,\,12, \\ 
 				\frac{3k^2+14k+16}{24} & \mbox{for} & k\equiv 10\,\,\mbox{mod}\,\,12, & 
 				\frac{3k^2+16k+13}{24} & \mbox{for} & k\equiv 11\,\,\mbox{mod}\,\,12.
			\end{array} 
			\right. 
		$$
	\end{satz}
	\textbf{Proof.}
	The list of graphs in Figure \ref{fig_graphs} is complete because the graphs have to be connected and the 
	complete graph on $4$ vertices is not a unit distance graph. Adding the formulas from Lemma 
	\ref{lemma_count_4_1} and Lemma \ref{lemma_count_4_2} yields the theorem.\hfill{$\square$}
	
\section{Computer enumeration of $\mathbf{k}$-polyominoes}
	
	For $n\ge 5$ we have constructed $k$-polyominoes with the aid of a computer and have obtained the 
	following values of $a_k(n)$ given in Table \ref{table_a_k_n} and Table \ref{table_a_k_n_2}. 
	
	\begin{table}[ht]
		\begin{center}
			\begin{tabular}[c]{|r|r|r|r|r|r|r|r|r|r|}
				\hline
				\n $k\!\!\setminus\!\! n$\n &  5\n & 6\n & 7\n & 8\n & 9\n & 10\n & 11\n & 12\n & 13\n \\
				\hline
				\n  3\n &\n   4\n &\n    12\n &\n     24\n &\n     66\n &\n     160\n &\n      448\n &	   
				\n    1186\n &\n      3334\n &\n     9235\n \\
				\n  4\n &\n  12\n &\n    35\n &\n    108\n &\n    369\n &\n    1285\n &\n     4655\n &
				\n  17073\n &\n      63600\n &\n    238591\n \\
				\n  5\n &\n  25\n &\n   118\n &\n    551\n &\n   2812\n &\n   14445\n &\n    76092\n &
				\n	 403976\n &\n   2167116\n &\n  11698961\n \\
				\n  6\n &\n  22\n &\n    82\n &\n    333\n &\n   1448\n &\n    6572\n &\n    30490\n &
				\n	 143552\n &\n    683101\n &\n   3274826\n \\
				\n  7\n &\n  25\n &\n   118\n &\n    558\n &\n   2876\n &\n   14982\n &\n    80075\n &
				\n  431889\n &\n   2354991\n &\n  12930257\n \\
				\n  8\n &\n  50\n &\n   269\n &\n   1605\n &\n  10102\n &\n   65323\n &\n	430302\n &
				\n  2868320\n &\n 19299334\n &\n 130807068\n \\
				\n  9\n &\n  82\n &\n   585\n &\n   4418\n &\n  34838\n &\n  280014\n &\n  2285047\n &
				\n 18838395\n &\n 156644526\n & \n 1311575691\n \\
				\n 10\n &\n 127\n &\n   985\n &\n   8350\n &\n  73675\n &\n  664411\n &\n  6078768\n & 
				\n 56198759\n &\n 523924389\n &       \\
				\n 11\n &\n 186\n &\n  1750\n &\n  17501\n &\n 181127\n &\n 1908239\n &\n 20376032\n &	
				\n 219770162\n &\n 2390025622\n &         \\		
				\n 12\n &\n 168\n &\n  1438\n &\n  13512\n &\n 131801\n &\n 1314914\n &\n 13303523\n & 
				\n 136035511\n &\n 1402844804\n &          \\
				\n 13\n &\n 187\n &\n  1765\n &\n  17775\n &\n 185297\n &\n 1968684\n &\n 21208739\n & 
				\n 230877323\n &           &           \\
				\n 14\n &\n 263\n &\n  2718\n &\n  30467\n &\n 352375\n &\n	4158216\n &\n 49734303\n &	
				\n 601094660\n &           &              \\
				\n 15\n &\n 362\n &\n  4336\n &\n  55264\n &\n 725869\n &\n 9707046\n &\n 131517548\n &	
				\n 1800038803\n &           &              \\
				\n 16\n &\n 472\n &\n  6040\n &\n  83252\n &\n 1180526\n &\n 17054708\n &\n 249598727\n &	
				\n 3690421289\n &           &            \\
				\n 17\n &\n 613\n &\n  8814\n &\n 134422\n &\n 2104485\n &\n 33522023\n &\n 540742895\n &
							  &           &             \\
				\n 18\n &\n 566\n &\n  7678\n &\n 112514\n &\n 1694978\n &\n 26019735\n &\n 404616118\n &
							  &           &             \\
				\n 19\n &\n 615\n &\n  8839\n &\n 135175\n &\n 2123088\n &\n 33942901\n &\n 549711709\n &
							  &           &            \\	
				\n 20\n &\n 776\n &\n 11876\n &\n 195122\n &\n 3291481\n &\n 56537856\n &\n 983715865\n &
							  &           &           \\
				\hline
			\end{tabular}
			\caption{Number of $k$-polyominoes with $n$ cells for small $k$ and $n$.}
			\label{table_a_k_n}
		\end{center}
	\end{table}
	
	\vspace*{-3mm}
		
	\noindent
	Now we go into more detail how the computer enumeration was done. At first we have to represent $k$-polyominoes 
	by a suitable data structure. As in Lemma \ref{lemma_overlapping} a $k$-polyomino can be described by the set of all 
	discrete angles between three neighbored cells. By fixing one direction we can define the discrete angle between this 
	direction and two neighbored cells and so describe a $k$-polyomino by an $n \times n$-matrix with integer 
	entries. Due to Corollary \ref{cor_neighbors} we can also describe it as a $6 \times n$-matrix by listing 
	only the neighbors. To deal with symmetry we define a canonical form for these matrices.
	
	Our general construction strategy is orderly generation \cite{winner}, where we use a variant introduced in 
	\cite{phd_kurz,paper_characteristic}. 
	Here a $k$-polyomino consisting of $n$ cells is constructed by glueing two $k$-polyominoes consisting of $n-1$ cells 
	having $n-2$ cells in 
	common. There are two advantages of this approach. In a $k$-polyomino each two cells must be nonoverlapping. If we would 
	add a cell in each generation step we would have to check $n-1$ pairs of cells whether they are nonoverlapping or not. By 
	glueing two $k$-polyominoes we only need two perform one such check. To demonstrate the the second advantage we compare in 
	Table \ref{table_compare} the numbers $c_1(n,k)$ and $c_2(n,k)$ of candidates produced by the original version and the 
	used variant via glueing of orderly generation.
	
	To avoid numerical twists in the overlapping check we utilize Gr\"obner bases \cite{groebner}. 
	
	\begin{table}[!ht]
		\begin{center}
			\begin{tabular}{|r|r|r|r|r||r|r|r|r|r|}
				\hline
				$\!k\!\!\setminus\!\! n\!$ & 5 & 6 & 7 & 8 & $\!k\!\!\setminus\!\! n\!$ &  5 & 	  6   & 7\\
				\hline
				21  &  972 &  16410 &  294091 &   5402087 & 36  &  4575 & 130711 &  3943836 \\
				22  & 1179 &  20970 &  397852 &   7739008 & 37  &  4796 & 140434 &  4326289 \\
				23  & 1437 &  27720 &  566007 &  11832175 & 38  &  5380 & 163027 &  5204536 \\
				24  & 1347 &  24998 &  495773 &  10079003 & 39  &  6089 & 193587 &  6464267 \\
				25  & 1439 &  27787 &  568602 &  11917261 & 40  &  6760 & 221521 &  7634297 \\
				26  & 1711 &  34763 &  751172 &  16624712 & 41  &  7578 & 259396 &  9311913 \\
				27  & 2045 &  44687 & 1031920 &  24389611 & 42  &  7282 & 244564 &  8643473 \\
				28  & 2376 &  54133 & 1307384 &  32317393 & 43  &  7584 & 259838 &  9341040 \\
				29  & 2786 &  67601 & 1729686 &  45260884 & 44  &  8373 & 295558 & 10958872 \\		
				30  & 2641 &  62252 & 1557663 &  39891448 & 45  &  9321 & 342841 & 13215115 \\
				31	 & 2790 &  67777 & 1737915 &  45587429 & 46  & 10207 & 385546 & 15274792 \\
				32  & 3204 &  81066 & 2169846 &  59424885 & 47  & 11282 & 442543 & 18169170 \\
				33  & 3706 &  99420 & 2808616 &  81124890 & 48  & 10890 & 420154 & 17012270 \\
				34  & 4193 & 116465 & 3413064 & 102292464 & 49  & 11290 & 443178 & 18217475 \\
				35  & 4789 & 140075 & 4306774 & 135337752 & 50  & 12309 & 495988 & 20944951 \\
				\hline
			\end{tabular}
			\caption{Number of $k$-polyominoes with $n$ cells for small $k$ and $n$.}
			\label{table_a_k_n_2}
		\end{center}
	\end{table}
	
	\begin{table}[!ht]
		\begin{center}
			\begin{tabular}[c]{|r|r|r|r|r|r|r|r|r|}
				\hline
				$n$ &  4 & 5 & 6 & 7 & 8 & 9 & 10 & 11 \\
				\hline
				\n $a_5(n)$\n    &\n   7\n &\n   25\n &\n   118\n &\n    551\n &\n    2812\n &\n    14445\n &\n    76092\n &\n 403976\n\\
				\n $c_1(n,5)$\n  &\n  21\n &\n   74\n &\n   242\n &\n   1038\n &\n    4476\n &\n    21945\n &\n   111232\n &\n 580139\n\\
				\n $c_2(n,5)$\n  &\n  19\n &\n   62\n &\n   192\n &\n    816\n &\n    3541\n &\n    17297\n &\n    87336\n &\n 452215\n\\
				\hline
				\n $a_7(n)$\n    &\n   7\n &\n   25\n &\n   118\n &\n    558\n &\n    2876\n &\n    14982\n &\n    80075\n &\n 431889\n\\
				\n $c_1(n,7)$\n  &\n  31\n &\n  107\n &\n   356\n &\n   1530\n &\n    6682\n &\n    33057\n &\n   168881\n &\n 889721\n\\
				\n $c_2(n,7)$\n  &\n  19\n &\n   62\n &\n   196\n &\n    821\n &\n    3584\n &\n    17778\n &\n    91109\n &\n 479814\n\\
				\hline
				\n $a_{13}(n)$\n &\n  23\n &\n  187\n &\n  1765\n &\n  17775\n &\n  185297\n &\n  1968684\n &\n 21208739\n &\n 230877323\n\\
				\n $c_1(n,13)$\n &\n 126\n &\n  721\n &\n  5059\n &\n  43842\n &\n  420958\n &\n  4294445\n &\n 45258582\n &\n 485481211\n\\
				\n $c_2(n,13)$\n &\n  76\n &\n  408\n &\n  2697\n &\n  23412\n &\n  223789\n &\n  2274489\n &\n 23849241\n &\n 254712159\n\\
				\hline
				\n $a_{17}(n)$\n &\n  48\n &\n  614\n &\n  8814\n &\n 134422\n &\n 2104485\n &\n 33522023\n &\n  540742895\n &\n       \n\\
				\n $c_1(n,17)$\n &\n 255\n &\n 2039\n &\n 22038\n &\n 292887\n &\n 4311681\n &\n 66600525\n &\n 1057440375\n &\n       \n\\
				\n $c_2(n,17)$\n &\n 171\n &\n 1261\n &\n 12964\n &\n 173839\n &\n 2545538\n &\n 39008006\n &\n  614066925\n &\n       \n\\
				\hline
			\end{tabular}
			\caption{Number of candidates $c_1(n,k)$ and $c_2(n,k)$ for $k$-polyominoes with $n$ cells.}
			\label{table_compare}
		\end{center}	
	\end{table}

	\section{Open problems for $\mathbf{k}$-polyominoes}
	For $4$-polyominoes the maximum area of the convex hull was considered in \cite{Bezdek1994}. If the area of 
	a cell is normalized to $1$ then the maximum area of a $4$-polyomino consisting of $n$ squares is given by 
	$n+\frac{1}{2}\left\lfloor\frac{n-1}{2}\right\rfloor\left\lfloor\frac{n}{2}\right\rfloor$. The second 
	author has proven an analogous result for the maximum content of the convex hull of a union of $d$-dimensional 
	units hypercubes \cite{dipl_kurz}, which is given by  
	$$
		\sum\limits_{I\subseteq\{1,\dots,d\}}\frac{1}{|I|!}\prod\limits_{i\in	I}\left\lfloor\frac{n-2+i}{d}\right\rfloor\
	$$
	for $n$ hypercubes. For other values of $k$ the question for the maximum area of the convex hull of $k$-polyominoes 
	is still open. Besides from \cite{0747.52010} no results are known for the question of the minimum area of the convex 
	hull, which is non trivial for $k\neq 3,4$. 	
	
	Another class of problems is the question for the minimum and the maximum number of edges of $k$-polyominoes. The following
	sharp inequalities for the number $q$ of edges of $k$-polyominoes consisting of $n$ cells were found in \cite{extremal} and
	are also given in \cite{0828.00001}.
	\begin{eqnarray*}
		k=3:&\quad& n+\left\lceil \frac{1}{2}\left(n+\sqrt{6n}\right)\right\rceil\le q\le 2n+1\\
		k=4:&\quad& 2n+\left\lceil 2\sqrt{n}\right\rceil\le q\le 3n+1\\
		k=6:&\quad& 3n-\left\lceil\sqrt{12n-3}\right\rceil\le q\le 5n+1
	\end{eqnarray*} 
	
	\noindent
	In general the maximum number of edges is given by $(k-1)n+1$. The numbers of $4$-polyominoes
	with a minimum number of edges were enumerated in \cite{counting}. 
	
	Since for $k\neq 3,4,6$ regular $k$-gons do not tile the plane the question about the maximum density $\delta(k)$ of an edge-to-edge 
	connected packing of regular $k$-gons arises. In \cite{phys}
	$$
		\delta(5)=\frac{3\sqrt{5}-5}{2}\approx 0.8541
	$$
	is conjectured.
	
	\bibliography{generalized_polyominoes}
	\bibdata{generalized_polyominoes}
	\bibliographystyle{plain}
\end{document}